\input amstex
\documentstyle {amsppt}
\tolerance=3000
\openup 6 pt
\nologo

\bigskip
\bigskip
\bigskip
\bigskip
\topmatter
\title
Forcing of periodic orbits for interval maps and\\
renormalization of piece wise affine maps.
\endtitle
\author
 Marco Martens \footnote {Institute of Mathematical Sciences, SUNY at 
Stony Brook, Stony Brook, NY 11794-3651.} 
and Charles Tresser\footnote {I.B.M., Po Box 218, Yorktown Heights, NY 10598.}
\endauthor
\endtopmatter
\bigskip
\bigskip
\bigskip
\bigskip
\bigskip
\bigskip
\centerline {\bf Abstract.}
\bigskip
 
\flushpar
We prove that for continuous maps on the interval, the existence
of a $n-$cycle, implies the existence of
$n-1$ points which interwind the original ones and are permuted by the map.
We then use this combinatorial result to show that piecewise affine maps 
(with no zero slope) cannot be infinitely renormalizable.

\newpage

\bigskip
\centerline{\bf 1. Introduction}
\bigskip

\flushpar
A fascinating feature of real analytic infinitely renormalizable
interval maps is that their attracting invariant Cantor sets seem to
have a complicated geometry (for the unimodal case, see
for instance $[S]$ and references therein). One could hope to avoid this
complexity by constructing piecewise affine examples (with no
zero slope). This is indeed
the case when there are infinitely many intervals of affinity (see e.g. 
$[T]$), 
but we show in section 3 that no example exists with finitely many intervals 
of affinity.
In order to prove this results we had to solve some questions about forcing
of permutations which are described in section 2. 

\flushpar
We now state our main results.
The collection of continuous maps on the interval is denoted by $C^0([0,1])$.
A {\it cycle } of a map $f\in C^0([0,1])$ is a collection
of pairwise disjoint closed intervals $\Cal I=\{I_0,I_1,..., I_q\}$ which are
cyclically permuted by $f$. That means $f(I_i)\subset I_{i+1 mod q}$. The 
cycle is called {\it trivial} if all intervals are points.

\flushpar
A cycle $\Cal{I}^j$ refines a cycle $\Cal{I}^i$
if $\cup I_l^j\subset \cup I_k^i$. In this situation there is a number
$a_{i,j}\in \Bbb{N}$ such that every component of $\Cal{I}^i$ contains 
$a_{i,j}$ components of $\Cal{I}^j$. We will always assume $a_{i,j}\ge 2$ and 
use the notation $\Cal{I}^i\supset\Cal{I}^j$ or $\Cal{I}^j\subset \Cal{I}^i$.

\flushpar
Let $\Cal{I}^j$ refine $\Cal{I}^i$. A connected component $G$ of 
$\cup I_l^i\setminus \cup I_k^j$ is called a {\it gap} if
$\partial G\subset \cup I_l^j$. The union of gaps is denoted
by $\Cal {G}(\Cal{I}^i,\Cal{I}^j)$.

\proclaim{Definition 1.1} Let $f\in C^0([0,1])$ and $\Cal{I}^2$ be a
cycle which refines the cycle $\Cal{I}^1$.

\flushpar
An invariant  set $P$ of periodic points of $f$ is called a {\it splitting }
of the pair of cycles
$\Cal{I}^1\supset\Cal{I}^2$ if

\parindent=15pt
\item{-} $P\subset \Cal {G}(\Cal{I}^1,\Cal{I}^2)$,
\item{-} every gap in $\Cal{G}(\Cal{I}^1,\Cal{I}^2)$ contains 
exactly one point of $P$.
\endproclaim

\proclaim{Theorem A} Every pair of cycles 
$\Cal{I}^1\supset \Cal{I}^2$ of a map $f\in C^0([0,1])$ admits 
a splitting.
\endproclaim

\flushpar
Theorem A is 
a main ingredient of the proof of Theorem C below.
It is a corollary of the following result.

\proclaim{Theorem B} If a continuous map $f$ on the interval has a periodic
orbit of period $n$ 
then it permutes $n-1$ points interwinding the periodic orbit.
\endproclaim

\flushpar
Here we say that a set of $n-1$ points on the real line {\it interwinds}
a set of $n$ points if any two consecutive points of any of these sets are
separated by a point of the other set. Theorem B in turns follows from 
a property of markov maps induced by permutations of points in the interval,
stated as Theorem 2.1.

\bigskip

\flushpar
A map $f\in C^0([0,1])$ is {\it infinitely renormalizable} if it has an 
infinite sequence of refining cycles
$$
\Cal{I}^1\supset\Cal{I}^2\supset \Cal{I}^3\supset\dots
$$
The set $\bigcap \cup \Cal{I}^n$ is called an {\it infinitely renormalizable
invariant set}.

\flushpar 
A map $f\in C^0([0,1])$ is called {\it piecewise affine} if there is a 
sequence of points $0=c_{-1}<c_0<\dots <c_{d}<c_{d+1}=1$, called 
{\it corner points}, such that $f$ is 
affine with non-zero slope on each interval 
$[c_i, c_{i+1}]$ and $f$ is not affine on any larger 
interval. The set of corners of $f$ is denoted by $C_f$.
The collection of all piecewise affine maps is denoted by $PL([0,1])$. 

\proclaim{Theorem C} There are no infinitely renormalizable piecewise
affine maps.
\endproclaim

\flushpar
Notice that the conclusion of Theorem C fails if $C_f$ is allowed to be
countable  (see e.g. $[T]$) or if we relax the non-zero slope condition
(see e.g. $[BMT]$).
Theorem C can be understood as a step toward proving the conjecture in $[GMT]$,
that piecewise affine maps on the interval without periodic attractors
are eventually expanding. An important ingredient of the proof of Theorem C
is the Expansion-Lemma in section 3. It states that there is a reasonable big 
collection of expanding periodic orbits with exponent away from $0$. This 
Expansion-Lemma should be compared with the Finiteness of Attractors Theorem 
in [MMS], which states that in smooth maps periodic orbits with sufficiently
high period, are expanding with  exponent away from $0$. Whether this is also
true for piecewise affine maps is part of the conjecture in $[GMT]$. 

\flushpar
{\bf Acknowledgements}. After proving Theorem C in the case of period
doubling ($a_{i,i+1}\equiv 2$) we checked with some colleagues about its 
originality. Some time later, Michal Misiurewicz and Karen Brooks reported
to us that V.J. Lopez and L. Snoha had recently obtained the same result
 [LS].

\tolerance=3000

\bigskip
\centerline{\bf 2. Forcing}
\bigskip

\flushpar
Let $S_n$ be the collection of permutations of $N_n=\{1,2,\dots,n\}$. For
every $\pi\in S_n$ we define the $n-1\times n-1$-matrix $F_\pi$, with
\parindent=15pt
\item{-} $F_\pi(i,j)=1$, if $(j,j+1)\subset (\pi(i),\pi(i+1))$,
\item{-} $F_\pi(i,j)=0$ otherwise.

\flushpar
A continuous map $f:[1,n]\to [1,n]$ with $f|N_n=\pi$ will in general
map the {\it gaps} $(i,i+1)$ in a very non-monotone way, and one cannot guess
the full set of periodic orbits of $f$ by the only knowledge of $\pi$. 
However $f|N_n=\pi$ implies some minimal complexity for the dynamics of $f$.
More specifically, it is known that the subshift of finite type defined by 
$F_\pi$ can always be monotonically imbedded into the dynamics of $f$;
we say that $\pi=f|N_n$ {\it forces} all the dynamics of the subshift
defined by $F_\pi$. In the above statement, monotonicity refers to the
skewed lexicographic order on the sequence space,
as used in kneading theory (see e.g.  [MT]),
 and the usual order on $[1,n]$.
As usual we identify the matrix $F_\pi$ with the corresponding subshift.

\flushpar
If $\phi\in S_m$ and $F$ is a $m\times m$-matrix with $0,1$ entries,
such that $F(i,\phi(i))=1$, we say that {\it $F$ contains the 
permutation $\phi$}, and write $\phi<<F$.

\proclaim{Theorem 2.1} For every $\pi\in S_n$, $n\ge 2$ there exists a 
$\phi\in S_{n-1}$ with
$$
\phi<<F_\pi.
$$
\endproclaim

\flushpar
\demo{Remark}
It is easy to construct examples of subshifts of finite types whose
defining matrices have some power with all entries positive but 
do not contain a permutation.
\enddemo

\flushpar
The proof of this Theorem needs some preparation. We are going to describe
a ``cutting'' procedure on permutations and a related ``cutting''
procedure on their matrices.
Let $n\ge 2$ and $k\le n-1$ and let 
$j_k:N_{n-1}\to N_n\setminus \{k\}$ stand for the order preserving bijection.

\flushpar
If  $\pi\in S_n$ and $k=\pi^{-1}(n)$ then $\hat{\pi}\in S_{n-1}$
is defined by
$$
\hat{\pi}(i)=\pi(j_k(i)),
$$
where $i\in N_{n-1}$.

\flushpar
From the definition of the matrix $F_\pi$, the $1$'s are consecutive in each
of its rows. Furthermore if $F_\pi(k,n)=1$ then $\pi(k)=n$ or
$\pi(k+1)=n$. So there are at most two rows with the last 
entry equal to $1$, and such rows have to be consecutive.

\flushpar
Consider the rows with last entry equal to $1$ and assume that row $k$ has the 
shortest block of $1$'s among those: that means, if $F_\pi(l,n)=1$ and $l\ne k$
then there exists $j\le n$
with $F_\pi(l,j)=1$ and $F_\pi(k,j)=0$.
Now we define a $(n-1)\times (n-1)$ matrix $\hat{F}_\pi$ by
$$
\hat{F}_\pi(i,j)=F_\pi(j_k(i),j).
$$
This $\hat{ }\,$-operation on matrices has been defined for matrices induced by
permutations. To the contrary of the $\hat{ }\,$-operation on permutations, 
in general it cannot be squared because the new matrix
is maybe not induced by a permutation (this would always be the case
for permutations realizable as restrictions of unimodal maps).
In particular $F_{\hat\pi}$ is in general not equal to $\hat{F}_\pi$. These
two matrices are only equal for unimodal permutations. In general we have

\proclaim{Proposition 2.2} For every $\pi\in S_n$, $n\ge 2$,
$$
F_{\hat\pi}\le \hat{F}_\pi.
$$
\endproclaim

\demo{Proof} Fix $\pi\in S_n$. For $k\le n-1$ let
$h_k:N_{n-2}\setminus \{k-1\}\to N_{n-1}\setminus \{k-1,k\}$ be
the order preserving bijection.

\flushpar
If $\pi^{-1}(n)=k$ then the $(k-1)^{th}$ and $k^{th}$ row have a $1$ in there
last entry. Denote these rows by the vectors $V$ and $v$. Assume that
$V$ has a longer block of $1's$. In the case that $k=1$ or $k=n$ there
is only one row whose last entry equals one, resp. the first or the 
last row. In these cases let $V$ be this row and $v=0$. 

\flushpar
Now the Proposition follows immediately from

\proclaim{Claim} If $i\ne k-1$ then
$$
\aligned
F_{\hat{\pi}}(i,j)=F_\pi(h_k(i),j),\\
\hat{F}_\pi(i,j)=F_\pi(h_k(i),j).
\endaligned
$$
For $j\le n-2$ and $k\ne 1,n$
$$
\aligned
F_{\hat{\pi}}(k-1,j)&=V(j)-v(j),\\
\hat{F}_\pi(k-1,j)&=V(j).
\endaligned
$$
\endproclaim

\demo{Proof}
The matrix $\hat{F}_\pi$ was obtained from $F_\pi$ by erasing the
last column and row $v$. The result of this operation is expressed
in the claim.

\flushpar
Now consider the $\hat{\pi}$ image of the gap $(i,i+1)$. If $i\le k-2$
then the $\hat\pi$-image is the same as the $\pi$-image. Hence the first
$k-2$ rows of $F_{\hat\pi}$ equal the first $k-2$ rows of $F_\pi$. If 
$i\ge k$ then the $\hat\pi$-image is the $\pi$-image of the gap 
$(i+1,i+2)$. So the last $(n-2)-(k-1)$ rows of $F_{\hat\pi}$ equal
the last $(n-1)-k$ rows of $F_\pi$. These properties are expressed
in 
$$
F_{\hat\pi}(i,j)=F_\pi(h_k(i),j),
$$
whenever $i\ne k-1$.
It just remains to consider the image of the $(k-1)^{th}$ gap. The boundary 
points are $k-1$ and $k$. Hence the $\hat\pi$-image of this gap equals
the interval $(\pi(k-1),\pi(k+1))$.
It follows easily that if a gap $(j,j+1)$ is covered  by the 
interval $(\pi(k-1),\pi(k+1)$ then $V(j)-v(j)=1$.
\hfill\hfill\qed $\,\,$ (Claim and Proposition 2.2)
\enddemo

\demo{Proof of Theorem 2.1}
The proof of Theorem 2.1 is by induction. For a permutation $\pi\in P_S$
we have that the only entry of $F_\pi$ equals $1$. Hence $F_\pi$ contains
a permutations. Now suppose that every matrix $F_\pi$ with $\pi\in S_n$,
contains a permutations.

\flushpar
Let $\pi\in S_{n+1}$. Then $\hat{\pi}\in S_n$ and hence $F_{\hat\pi}$ contains a 
permutation $\phi_0<<F_{\hat\pi}$. From Proposition 2.2 we get 
$F_{\hat\pi}\le \hat{F}_\pi$. Hence $\phi_0<<\hat{F}_\pi$.

\flushpar
Say, that $\hat{F}_\pi$ was obtained by cutting the $k^{th}$ row of 
$F_\pi$. Then $F_\pi$ contains the permutation $\phi\in S_{n+1}$
defined by
$$
\aligned
\phi(k)&=n \\
\phi(i)&=\phi_0(j_k^{-1}(i)) \text{ if } i\ne k.
\endaligned
$$
\hfill\hfill\qed $\,\,$ (Theorem 2.1)
\enddemo

\demo{Proof of Theorem B} Let $f$ be a continuous map with a periodic orbit
$\{p_1,p_2,\dots,p_n\}$ of period $n$. Assume $p_1<p_2<\dots<p_n$. And let
$G_i=(p_i,p_{i+1})$, with $i=1,2,\dots, n-1$.
By Theorem 2.1 we know that there exists a permutation $\phi \in S_{n-1}$
and intervals $T_i\subset G_i$ such that $f:T_i\to G_{\phi(i)}$ is onto.
The continuity of $f$ assures the existence of points $s_i\in \overline{T_i}$ 
such that $f(s_i)=s_{\phi(i)}$, with $i=1,2,\dots,n-1$. Clearly the points 
$s_i$ are periodic with period less than $n$. In particular 
$s_i\in int(G_i)$, they actually interwind the original orbit.

\hfill\hfill\qed $\,\,$ (Theorem B)
\enddemo

\demo{Proof of Theorem A} To prove Theorem A we may collapse the cycle 
$\Cal{I}^2$ into a periodic orbit, say $\Cal{I}^2= \{p_1,p_2,\dots,p_{q_2}\}$.
Again,
as in the proof of Theorem B, let $G_i=(p_i,p_{i+1})$. 
Apply Theorem 2.1 to get a permutation
$\phi\in S_{q_2-1}$ and the intervals $T_i\subset G_i$ such that
$f:T_i\to G_{\phi(i)}$ is onto. Let $P'=\{s_1,s_2,\dots, s_{{q_2}-1}\}$ be
the corresponding points which interwind $\{p_1,p_2,\dots, p_{q_2}\}$, 
that means  $s_i\in T_i$ and $f(s_i)=s_{\phi(i)}$.

\flushpar
Unfortunately $P'$ will not be a splitting for the pair 
$\Cal{I}^2\subset \Cal{I}^1$, it contains also points outside $\Cal{I}^2$.
Let $P=P'\cap G(\Cal{I}^1,\Cal{I}^2)$. To show that $P$ is a splitting
it suffices to show that $P$ is invariant. Take $s_i\in P$. Then 
$T_i$ is subset of the cycle $\Cal{I}^1$. Hence $f(T_i)=G_{\phi(i)}$
is subset of the cycle $\Cal{I}^1$, in particular $G_{\phi(i)}$ is a gap. 
So $s_{\phi(i)}\in G(\Cal{I}^1,\Cal{I}^2)\cap P'=P$.

\hfill\hfill\qed $\,\,$ (Theorem A)
\enddemo

\tolerance=3000

\bigskip
\centerline{\bf 3. Renormalization}
\bigskip

\flushpar
In his section we are going to prove Theorem C. The main reason why
piecewise affine maps cannot be infinitely renormalizable, is the fact
that there are enough periodic orbits with some definite expansion.

\flushpar
A cycle $\Cal{I}^2$ is called a {\it doubling} of a cycle $\Cal{I}^1$ if
$\Cal{I}^1\supset \Cal{I}^2$ and $a_{1,2}=2$.

\proclaim{Expansion-Lemma 3.1} Let $f\in PL([0,1])$ having two cycles
$\Cal{I}^1\supset \Cal{I}^2$ with splitting $P$. If 
\parindent=15pt
\item{-} $\Cal{I}^2$ is not a doubling of some $\Cal{I}\subset\Cal{I}^1$ and
\item{-} $C_f\subset \cup\Cal{I}^2$ 

then there exists $x\in P$ with
$$
|Df^p(x)|\ge 1+e^{-11V},
$$
where $p$ is the period of $x$ and V the variation of $\log(|Df|)$.
\endproclaim

\demo{proof} Let $I\in \Cal{I}^1$ and $G\subset I$ a gap, say 
$\partial G \subset I^2_l\cup I^2_r$. Now there exists
a periodic point $x\in  P\cap G$. Say, it has period $p\in \Bbb{N}$.

\proclaim{Claim 1} There exists a gap $G'\subset I$ and $G'\ne G$ and an 
interval $K\subset G$ with $f^p: K\to G'$ onto. 
\endproclaim

\demo{Proof Claim 1}
Say, $\Cal{I}^2$ has $q_2$ intervals. So $p\le q_2-1$. 
Assume Claim 1 is not true then
$f^p(G)\subset I^2_l\cup G\cup I^2_r$ and 
$f^p(I^2_l\cup I^2_r)\subset I^2_l\cup I^2_r$. Because $p<q_2$,
the two intervals $I^2_l$ and $I^2_r$ have to be interchanged. In particular
$p=\frac12 q_2$. 
Let $H=I^2_l\cup G\cup I^2_r$ then $f^p(H)\subset H$. If we can show that
the orbit $\{H,f(H),\dots,f^{p-1}(H)\}$ is pairwise disjoint, then we would
have shown that $\Cal{I}^2$ is a doubling of some cycle 
$\Cal{H}\subset \Cal{I}^1$. This contradicts the assumptions and Claim 1
would be proved.

\flushpar
So suppose $f^j(H)\cap f^p(H)\ne \emptyset$, for some positive $j\le p-1$. 
Now $I^2_{l+j}$
and $I^2_{r+j}$ do not intersect $H$, since only $I^2_l$ and $I^2_r$ intersect
$H$. Assume with no loss of generality $I^2_l\subset f^j(H)$. This implies 
$I^2_{l+p-j}=f^{p-j}(I^2_l)\subset 
f^p(H)\subset H$. 
This is impossible because only $I^2_l$
and $I^2_r$ are in $H$. We proved that the orbit of $H$ forms a cycle.
\hfill\hfill\qed $\,\,$ (Claim 1)
\enddemo

\demo{Reamrk} In Claim 1, $G'$ can be considered to be a gap adjacent to 
$I^2_l$ or $I^2_r$.
\enddemo

\flushpar
Now let $T\subset G$ be the maximal interval containing $x$ such that
$(f^p(T)\cap \cup \Cal{I}^2)\subset I^2_l\cup I^2_r$. Using Claim 1
we find an interval $K\subset T$ with $f^p:K\to G'$ onto. 
Because all corner points and their orbits
are in $\cup \Cal{I}^2$, the map $f^p:K\to G'$ is in 
fact affine and onto. A collection of intervals in $[0,1]$ is said to have 
{\it intersection multiplicity $w$} if every point in $[0,1]$ is contained
in at most $w$ intervals of the collection.

\proclaim{Claim 2} The intersection multiplicity of
$\{T,f(T),f^2(T),\dots,f^{p-1}(T)\}$ is at most $11$.
In particular 
$$
Var(\log(|Df^p|T|))\le 11 V.
$$
\endproclaim

\demo{Proof of Claim 2} To prove Claim 2 it is enough to show that there
are at most $10$ values $1\le i\le p-1$ such that $f^i(T)\cap f^p(T)\ne\emptyset$.

\flushpar
Let $I^2_L$ be the left neighbor of $I^2_l$ and $I^2_R$ be the right neighbor 
of $I^2_r$. The interval between $I^2_L$ and $I^2_R$ is denoted by $S$. Clearly
$f^p(T)\subset S$. Observe that $S$ contains at most $3$ gaps.
Because the orbit of $x$ is one of the splitting periodic orbits in $P$, 
the orbit of $x$ intersects $S$ in at most $3$ points, say in $x$, 
$f^a(x)$ and $f^b(x)$. These three intersections can also give rise to an 
intersection of $f^a(T)$ or $f^b(T)$ with $f^p(T)$.

\flushpar
Consider an intersection $f^i(T)\cap f^p(T)\ne \emptyset$, where  
$i\le p-1$ and
$i\notin\{0,a,b\}$. Because $f^i(x)\notin I^2_L\cup S\cup I^2_R$,
the interval $I^2_L$ (or $I^2_R$) is contained in $f^i(T)$. 
Hence $L+p-i\in \{l,r,l+q_2,r+q_2\}$
(or $R+p-i\in \{l,r,l+q_2,r+q_2\}$). 
This means that there are at most $8$ possible
values for $i\notin \{0,a,b\}$ giving rise to an intersection. 

\flushpar
All together we get at most $2+4+4=10$ intersections 
$f^i(T)\cap f^p(T)\ne\emptyset$, with $1\le i\le p-1$. The intersection 
multiplicity is at most $11$.
\hfill\hfill\qed $\,\,$ (Claim 2)
\enddemo

\flushpar
Because the orbits of the corners are contained in $\cup \Cal{I}^2$
there exists an interval $D\subset T$ with $x\in D$ and $f^p:D\to G$
affine and onto.

\flushpar
To prove Lemma 3.1, assume that $G$ was chosen the smallest
gap in $I$, $|G|\le |G'|$. Furthermore $|K|\le |T|-|D|$.
Observe $\frac{|G|}{|D|}=|Df^p(x)|$ and 
$\frac{|G'|}{|K|}=|Df^p_{|K}|$.
Then
$$
\aligned
11V&\ge \log(|Df^p|K|)-\log(|Df^p(x)|)\\
  &=\log(\frac{|G'|}{|K|}\times\frac{|D|}{|G|})\\
  &\ge \log\frac{|D|}{|K|}\\
  &\ge \log\frac{|D|}{|T|-|D|}.
\endaligned
$$
This implies
$$
|Df^p(x)|=\frac{|G|}{|D|}\ge \frac{|T|}{|D|}\ge 1+e^{-11 V}.
$$
\hfill\hfill\qed $\,\,$ (Lemma 3.1)
\enddemo

\demo{Proof of Theorem C}
The proof will be given in $PL$. This
is the collection of piecewise affine maps, not defined on the interval,
but defined on a finite union of intervals. The proof of Theorem C will
be by contradiction.  
Suppose that $PL([0,1])$ has an infinitely renormalizable map. Then also
$PL$ has an infinitely renormalizable map. Let $f\in PL$ be an
infinitely renormalizable map, whose number of corner points is minimal.
Denote the cycles by
$$
\Cal{I}^1\supset\Cal{I}^2\supset \Cal{I}^3\supset\dots
$$ 
Use the notation $q_n=\#\Cal{I}^n$.
We may assume that this sequence of cycles is complete. This means that
if there is some cycle $\Cal{I}$ with $\Cal{I}^{n+1}
\subset\Cal{I}\subset \Cal{I}^n$ then $\#\Cal{I}=q_n$ or
$\#\Cal{I}=q_{n+1}$. It can be shown that the infinitely renormalizable
invariant set is a minimal Cantor set (this follows from the non-existence
of wandering intervals for maps in $PL$. The minimality of the action of $f$
on its infinitely renormalizable Cantor set and the minimality of 
the number of corner points allow us to assume 

\proclaim{Claim 1} $C_f\subset int(\cup\Cal{I}^n)$ for all 
$n\in\Bbb{N}$. \hfill\hfill\qed $\,\,$ (Claim 1)
\endproclaim

\flushpar
Let $V=Var(\log(|Df|))$ and set 
$R_{n,j}=f^{q_n}|I^n_j$, $n\in\Bbb{N}$ and $j=0,1,2,\dots,q_n-1$.

\proclaim{Claim 2} There exists $K\in\Bbb{R}$ such that
for all $n\in \Bbb{N}$, $j<q_n$
$$
|DR_{n,j}(x)|\le K,
$$
for all $x\in I^n_j$. 
\endproclaim

\demo{proof}
Observe that $R_{n,j}\in PL(I^n_j)$ and that it has uniform, that is
independent of $n$ and $j$, bounds on the numbers of corner points.
Furthermore
$$
Var(\log(|DR_{n,j}|))\le V.
$$
Now Claim 2 follows easily.
\hfill\hfill\qed $\,\,$ (Claim 2)
\enddemo

\flushpar
Let $P_n$ be a splitting for the pair $\Cal{I}^n\supset\Cal{I}^{n+1}$. 
Claim 1 allows us to define
$$
\aligned
B_n&=\sum_{x\in \partial (\cup\Cal{I}^n)} \log(|Df(x)|),\\
M_n&=\sum_{x\in P_n}\log(|Df(x)|).
\endaligned
$$

\proclaim{Claim 3} For every $n\in\Bbb{N}$
$$
B_{n+1}=B_n+2M_n.
$$
In particular
$$
B_n=B_1+\sum_{k=1}^{n-1} 2M_n. 
$$
\endproclaim

\demo{Proof} Consider $y\in \partial I^{n+1}_j\subset I^n_i$. 
If $y$ is also
in the boundary of some gap $G$, then $Df(y)=Df(x)$, where 
$x\in P_n\cap G$. This equality holds because all corner points are
in $\cup\Cal{I}^{n+1}$. In the other case, $y$ is not in the boundary of 
some gap then $Df(y)=Df(z)$, where $z\in \partial I^n_i$. Again equality
holds because $C_f\subset \cup\Cal{I}^{n+1}$.

\flushpar
Observe that all boundary points of $\cup\Cal{I}^n$ are counted once and all
splitting periodic points twice. 
\hfill\hfill\qed $\,\,$ (Claim 3)
\enddemo

\proclaim{Claim 4} For every periodic point $x$, with period $p$,
$\log(|Df^p(x)|>0$. In particular 
$$
M_n>0,
$$
for all $n\in \Bbb{N}$.
\endproclaim

\demo{proof} It is easy to see that a periodic attractor attracts
a corner point. Because $q_n\ge 2^{n-1}q_1\to \infty$ every periodic 
orbit falls eventually outside $\cup\Cal{I}^n$. An attracting periodic
orbit would take a corner point with it, contradicting the minimality
of corner points of $f$. 
The non-existence of neutral periodic orbits follows by a similar
argument.
\hfill\hfill\qed $\,\,$ (Claim 4)
\enddemo

\flushpar
By Claims 3 and 4 we know that the sequence $B_n$ increases.

\proclaim{Claim 5}
$$
\lim_{n\to\infty} B_n=B<\infty.
$$
\endproclaim

\demo{Proof} The corners of $f$ are denoted by $C_f=\{c_1,c_2,\dots,c_d\}$.
Take $x\in I^n_0$. Furthermore let $c_j\in I^n_{k_j}$ 
and $v_j=\log(|Df_{+}(c_j)|)-\log(|Df_{-}(c_j)|)$, $j=1,2,\dots,d$ (the $+$
and $-$ means right and left derivatives).
From [MMS] we know that $f$ does not have wandering intervals. This implies
that the length of the intervals in $\Cal{I}^n$ tends to zero for $n\to\infty$.
So (in $PL$) 
we may assume that every interval in $\Cal{I}^n$ contains at most one
corner point.

\flushpar
It is an easy exercise to compare $\log(|Df^{q_n}(x)|)$ with $B_n$. There
exist numbers $\sigma_{n,j}(x)\in \{-1,1\}$ such that
$$
2 \log(|Df^{q_n}(x)|)=  B_n+ \sum_{j=1}^d \sigma_{n,j}(x) v_j,
$$
where $\sigma_{n,j}(x)=1$ if and only if $f^{k_j}(x)\in I^2_{k_j}$ is on 
the right of $c_j$ (remember $c_j$ is the unique corner point in $I^2_{k_j}$).

\flushpar
From Claim 2 we get a bound on $|Df^{q_n}(x)|$. The finite sum in the
right hand side of the above equation is clearly uniformly bounded. Hence
we get a uniform bound on $B_n$.

\hfill\hfill\qed $\,\,$ (Claim 5)
\enddemo

\flushpar
As a consequence we get
$$
\lim_{n\to\infty } M_n=0.
$$
But now we can apply the Expansion-Lemma: eventually $\Cal{I}^{n+1}$
is a doubling of $\Cal{I}^n$, for all $n\ge n_0$. This means that the 
splitting $P_n$ becomes a single periodic orbit hitting every component
of the cycle $\Cal{I}^n$ exactly once.

\flushpar
As above we get the existence of numbers $\sigma_{n,j}\in\{-1,1\}$ such that
$$
2M_n=B_n+\sum_{j=1}^d \sigma^n_j v_j.
$$
Because $M_n\to 0$ and $B_n\to B$ and the fact that the sum in the above
equality takes only finitely many values we get eventually
$$
\sum_{j=1}^d \sigma^n_j v_j= -B.
$$
Combining this with $B_{n+1}=B_n+2M_n$ we get
$$
B_{n+1}=2B_n-B.
$$
For the sequence $B_n$ to be bounded we need $B_n\equiv B$.
So $M_n\equiv 0$, contradicting Claim 4.

\hfill\hfill\qed $\,\,$ (Theorem C)
\enddemo

\tolerance=3000

\bigskip
\centerline{\bf References}
\bigskip

\parindent=40pt
\item{[BMT]} K.M. Brucks, M. Misiurewicz, C. Tresser, {\it Monotonicity
             properties of the family of trapezoidal maps}, Commun. Math.
             Phys. {\bf 137} (1991) 1-12.
\item{[GMT]} R. Galeeva, M. Martens, C. Tresser, {\it Inducing, Slopes, and
             Conjugacy Classes}, preprint 1994/4 at SUNY at Stony Brook.
\item{[MMS]} M. Martens, W. de Melo, S. van Strien, {\it Julia-Fatou-Sullivan
             Theory for real one-dimensional dynamics}, Acta Math. {\bf 168 }
             (1992) 273-318.
\item{[MT]} J. Milnor, W. Thurston, {\it On iterated maps of the interval},
            Springer Lecture Notes in Mathematics {\bf 1342} (1988) 465-563.
\item{[LS]} V.J. Lopez, L. Snoha, to appear.
\item{[S]}  D. Sullivan, {\it Bounds, quadratic differentials and renormalization
            conjectures}, Mathematics into the twenty-first Century: 1988
            Centenial Symposium, ed. F. Browder, Amer. Math. Soc. (1992) 417-466.
\item{[T]}  C. Tresser, {\it Fine structure of universal Cantor sets},
            Instabilities and Nonequilibrium structures III, eds. E. Tirapegui and
            W. Zeller (Reidel) Dordrecht (1991).

\bye